\documentclass[11pt]{article}
\usepackage{pst-node,pstricks,pst-plot,pst-text}%
\usepackage{latexsym,amssymb}

\font\tensym=msbm10
\font\sevensym=msbm7
\font\fivesym=msbm5

\font\tengoth=eufb10
\font\sevengoth=eufb7
\font\fivegoth=eufb5

\newfam\symfam
\textfont\symfam=\tensym
\scriptfont\symfam=\sevensym
\scriptscriptfont\symfam=\fivesym
\def\sym{\fam\symfam\tensym}

\newfam\gothfam
\textfont\gothfam=\tengoth
\scriptfont\gothfam=\sevengoth
\scriptscriptfont\gothfam=\fivegoth
\def\goth{\fam\gothfam\tengoth}

 \newtheorem{thm}{Theorem}
\newtheorem{prop}[thm]{Proposition}
\newtheorem{lem}[thm]{Lemma}
\newtheorem{cor}[thm]{Corollary}
\newtheorem{defn}[thm]{Definition}
\newtheorem{note}{Note}

\def\cqfd{\hfill \\ \hbox{ }\hfill$\Box$ \bigskip}

\def\bb{\sym}
\def\B{{\sym B}}

\def\Z{{\sym Z}}
\def\M{{\sym M}}
\def\F{{\sym F}}

\def\GG{{\sym G}}
\def\L{{\sym L}}
\def\OV#1{\overline{#1}}

\def\R{{\bb R}}
\def\N{{\bb N}}

\def\Cont{\mbox{\rm Cont }}
\def\CQFD{\\ \hbox{ }\hfill$\Box$}
\def\preceq{\underline{\prec}}
\def\Alph{\mbox{Alph}}
\newtheorem{remark}{Remark}{\itshape}{\rmfamily}
\newtheorem{example}{Example}{\itshape}{\rmfamily}

\begin{document}
\bibliographystyle{plain}
\title{\bf Transitive factorizations of free partially
commutative monoids and Lie algebras}
\author{G\'erard Duchamp and Jean-Gabriel Luque\\
LIFAR, Facult\'e des Sciences et des Techniques,\\ 76821
Mont-Saint-Aignan CEDEX, France.}

\maketitle
\begin{abstract}
Let $\M(A,\theta)$ be a free partially commutative monoid. We give here a
necessary and sufficient condition on a subalphabet $B\subset A$ such that
the right factor of a bisection $\M(A,\theta)=\M(B,\theta_B).T$ be also
partially commutative free. This extends strictly the (classical)
elimination theory on partial commutations and allows to construct new
factorizations of $\M(A,\theta)$ and associated bases of $L_K(A,\theta)$.
\end{abstract}
\begin{center}{\small\bf R\'esum\'e}\end{center}
\begin{quote}
\small Soit $\M(A,\theta)$ un mono\"\i de partiellement commutatif libre.
Nous donnons une condition n\'ecessaire et suffisante sur un sous alphabet
$B\subset A$ pour que le facteur droit d'une bisection de la forme
$\M(A,\theta)=\M(B,\theta_B).T$ soit partiellement commutatif libre. Ceci
nous permet d'\'etendre strictement  et de fa\c con optimale  la th\'eorie
(classique) de
l'\'elimination avec commutations partielles et de construire de nouvelles
factorisations de $\M(A,\theta)$ ainsi que les bases de $L_K(A,\theta)$
associ\'ees.
\end{quote}
\baselineskip = 15pt
\section{Introduction}

A factorization of a monoid is a direct decomposition
$$
M=\prod^{\leftarrow}_{i\in I} M_i
$$
where $M$ and the $M_i$ are monoids and $I$ is totally ordered. This notion
is due to {\sc Sch\"utzenberger} (see \cite{SC1,Sc} where the link with the
free Lie algebra is studied). Then, in his Ph. D. \cite{Vi1}, {\sc Viennot}
showed how combinatorial bases of the free Lie algebra could be constructed
by composition of bisections (i.e. $|I|=2$) obtained by elimination of
generators (ideas initiated by {\sc Lazard} \cite{Laz} and {\sc Shirshov}
\cite{Shi}). One of the authors with D. Krob found similar decompositions
for the free partially commutative monoid into free factors and studied the
link with Lie algebras and groups \cite{DK1}. This works generalizes the
completely free case, but  is 
restricted to the situation where the outgoing factors are also free. \\
Here, we study the general problem of eliminating generators in these
structures and first remark that in any (set theoretical) direct
decomposition
$$
M(A,\theta)=M(B,\theta_B).T
$$
(with $B\subset A$, a subalphabet) the complement is a monoid. We get a
criterion
to characterize the case when $T$ is free partially commutative and
construct bases of the associated Lie algebras. The case of the group is
also mentionned.

\section{Definitions and background}
We recall that the free partially commutative monoid is defined by
generators and relations as
$$
\M(A,\theta)=\langle A|ab=ba, (a,b)\in \theta\rangle_{Mon},
$$
where $A$ is an alphabet and $\theta\subset A\times A$ is an antireflexive
(i.e. without loops) and symmetric graph on $A$ ($\theta$ is called an
independence relation). Thus, $\M(A,\theta)$ is the quotient
$A^*/_{\equiv_\theta}$ where $\equiv_\theta$ is the congruence generated by
the set $\{(ab,ba)|(a,b)\in \theta\}$.
\begin{defn}\label{ComThe}
If $X$ is a subset of $\M(A,\theta)$, we set
$$\theta_{X}=\{(x_1,x_2)\in X^2 | \Alph(x_1)\times \Alph(x_2)\subset
\theta\}.$$
\end{defn}
Note that $(x_1,x_2)\in\theta_X$ implies $\Alph(x_1)\cap 
\Alph(x_2)=\emptyset$, similarly
we
set $\theta_\M=\theta_{\M(A,\theta)}$.\\
As in \cite{DK2}, we denote $IA(t)=\{z\in A|t=zw\}$ and $TA(t)=\{z\in
A|t=wz\}$. \\
If $X$ is a subset of $\M(A,\theta)$, we denote $\langle X\rangle$ the
submonoid generated by 
$X$. In \cite{CH} and \cite{DR}, Choffrut introduces the partially
commutative
codes as some generating sets of free partially commutative submonoids. Let
$X$
be a set, we can prove easily that this definition is equivalent to the
fact that each trace $t\in\langle X\rangle$ admits a unique decomposition
on $X$ up to the commutations (i.e. $(X,\theta_X)$ is the independence
alphabet of $\langle X\rangle$ the submonoid generated by $X$).
\begin{example}
(i) Each subalphabet $B$ of $A$ is a partially commutative code.\\
(ii) Let $(A,\theta)= a-b\quad c$. The set $\{c,cb,ca\}$ is a code but not
the set $\{b,a,ca,cb\}$.
\end{example}
\section{Transitive bisections}
\subsection{Generalities}
We recall the definition of a factorization in the sense of
Sch\"utzenberger (cf. Viennot in \cite{Vi1} and \cite{Vi2}), this notion
will be reused extensively at the end of the paper.
\begin{defn}
(i) Let $\M$ be a monoid and $(\M_i)_{i\in J}$ an ordered family of
submonoids (the total ordering on $J$ will be denoted $<$). The family
$(M_i)_{i\in J}$ will be called a {\rm factorization} of $\M$ if and only
if every $m\in\M^+=\M-\{1\}$ has a unique decomposition
\[m=m_{i_1}m_{i_2} \dots m_{i_n}\]
with $i_1>i_2> \dots >i_k$ and for each $k\in[1..n]$,
$m_{i_k}\in\M^+_{i_k}$.\\
(ii) In the case of a free partially commutative monoid, a factorization
will be denoted by the sequence of the minimal generating sets of its
components.
\end{defn}
In the maximal case (each monoid has a unique generator), the factorization
is called {\it complete}.
\begin{example} (Complete factorizations in free and free partially
commutative
monoids.)\\
In the free monoid, it exists many complete factorizations. The most famous
of this kind
being the Lyndon factorization (defined as the set of primitive words
minimal in their conjugacy classes) is an example of a complete
factorization \cite{Lo,Re,SC1}. Hall sets defined in \cite{Sc} give us a
wider
example.\\
The set of Lyndon traces (i.e. the generalization of Lyndon words to the
partially commutative case, defined by Lalonde in \cite{La}) endowed with
the lexicographic ordering is a complete factorization of the free
partially
commutative monoid.
\end{example}
In the smallest case ($|J|=2$), the factorization is called a {\it
bisection}. Let $M$ be a monoid, then $(M_1,M_2)$ is a bisection of $M$ if
and only if
the mapping
\[M_1\times M_2\rightarrow M
\]
\[(m_1,m_2)\rightarrow m_1m_2\]
is one to one.\\
\begin{remark}
Not every submonoid is a left (right) factor of a bisection. If $M=M_1M_2$
is a bisection then $M_1$ satisfies $(u,uv\in M_1)\Rightarrow (v\in M_1)$
(see \cite{Dub}), however, this condition is not sufficient as shown by
$M_1=2\Z\subset\Z=M$.
\end{remark}
In case $M=\M(A,\theta)$, one can prove the following property.
\begin{prop}\label{PFACT2}
Let $(A,\theta)$ be an independence relation and $B\subset A$. Then
$\M(B,\theta_B)$ is the left (resp. right) factor of a bisection of
$\M(A,\theta)$.
\end{prop}
{\bf Proof } We treat, here, only the left case, the right case being
symmetrical.\\
It is clear that $M=\{t\in\M(A,\theta)|IA(t)\subset A-B\}$ is always a
monoid and that we have 
the (set theoretical) equality $\M(A,\theta)=\M(B,\theta_B).M$. It suffices
to prove the unicity of the 
decomposition of a trace. Let $w,w'\in\M(B,\theta_B)$ and $t,t'\in M$ such
that $wz=w'z'$. Using 
Levi's lemma, we find four traces $p,q,r,s$ such that $w=ps$, $t=rq$,
$w'=pr$ and $t'=sq$. But, 
by definition of $M$, we have $uv\in M$ implies $u\in M$, then $r,s\in 
M\cap\M(B,\theta_B)=\{1\}$. It follows $w=w'$ and $t=t'$, which gives the
result. \CQFD\\
In the sequel, we denote $Z=A-B$.\\
In the left case, the right submonoid above has
\[\beta_Z(B)=\{zw/z\in Z, w\in \M(B,\theta_B), IA(zw)= \{z\}\}\]
as  minimal generating subset.
\begin{remark}The monoid $\langle \beta_Z(B)\rangle$ may not be free
partially commutative. For example, if $A=\{a,b,c\}$,
\[ \theta: a-b\quad c \]
and $B=\{c\}$ then $a,b,ac,bc\in \beta_Z(B)$ and $a.bc=b.ac$.\\ \\
\end{remark}
\subsection{Transitively factorizing subalphabet}
Here we discuss a criterium for the complement $\langle \beta_Z(B)\rangle$
to be a free partially commutative submonoid.
\begin{defn}
Let $B\subset A$, we say that $B$ is a {\rm transitively factorizing
subalphabet} (TFSA) if and only $\beta_Z(B)$ is a partially commutative
code.
\end{defn}
We prove the following theorem.
\begin{thm}
Let $B\subset A$.
The following assertions are equivalent.
\begin{enumerate}
\item[(i)]The subalphabet $B$ is a TFSA.
\item[(ii)]The subalphabet $B$ satisfies the following condition.\\
For each $z_1\neq z_2\in Z$ and $w_1,
w_2, w'_1, w'_2\in\M(A,\theta)$ such that $IA(z_1w_1)=IA(z_1w'_1)=\{z_1\}$
and
$IA(z_1w_2)=IA(z_2w'_2)=\{z_2\}$ we have
\[z_1w_1z_2w_2=z_2w'_2z_1w'_1\Rightarrow w_1=w'_1, w_2=w'_2. \]
\item[(iii)] For each $(z,z')\in Z^2\cap \theta$, the
dependence\footnote{The
dependence graph is defined by $A\times A-\Delta-\theta$ where
$\Delta=\{(a,a)/a\in A\}$.} graph has no partial graph\footnote{We repeat
here the notion of partial graph. A graph $G'=(S',A')$ is a partial graph
of $G=(S,A)$ if and only if $S'\subset S$ and $A'\subset A\cap S'\times S'$
($G'$ is a subgraph of $G$ when equality $S=S'$ occurs).} like
\[z-b_1-\dots-b_n-z'.\]
with $b_1,\dots,b_n\in B$.
\end{enumerate}
\end{thm}
{\bf Proof}
It is easy to see that (i)$\Rightarrow$(ii) : by contraposition, if $B$
does not satisfy (ii) we can find $z_1w_1,z_2w_2, z_1w'_1, z_2w'_2\in
\beta_Z(B)$ such
that $z_1w_1.z_2w_2=z_2w'_2.z_1w'_1$ with $w_1\neq w'_1$ or $w_2\neq w'_2$
and this implies obviously that $\beta_Z(B)$ is not a partially commutative
code.\\
Let us prove that (ii)$\Rightarrow$(iii). Suppose that
\[z-b_1-\dots-b-n-z'\]
is a partial graph of the dependence graph, then it exists a subgraph of the
dependence graph of the form
\[z-c_1-\dots-c_m-z'\]
with $c_i\in B$. Consider the smallest integer $k$ such that
$(c_{k+1},z')\not\in\theta$. Then we have $zc_1\dots c_k.z'c_{k+1}\dots
c_{m}=z'.zc_1\dots c_m$, which proves that $B$ does not satisfy (ii).\\

Finally, we prove that (iii)$\Rightarrow$ (i). For each $z\in Z$, we define
$B_z$ the set of 
letters of $B$ having in the dependence graph a path leading to $z$ and all
inner points 
belonging to $B$. Clearly the assertion (iii) is equivalent to the fact that
$(z,z')\in \theta_Z$ implies 
$(\{z\}\cup B_z)\times (\{z'\}\cup B_{z'})\subseteq\theta$. It follows that
$\beta_z(B)\times\beta_{z'}(B)\subset\theta_\M$.\\
Consider the mapping $\mu$
from $Z$ into $K\langle\langle A,\theta\rangle\rangle$ defined by
$\mu(z)=\underline{\beta_z(B)}$. As $(z,z')\in \theta_Z\Rightarrow
[\mu(z),\mu(z')]=[\underline{\beta_z(B)},\underline{\beta_{z'}(B)}]=0$ and
$\langle \mu(z),1\rangle=0$\footnote{Here, for a series $S=\sum\alpha_uu$,
we denotes $\langle 
S,w\rangle=\alpha_w$.}, we can extend $\mu$ as a continuous morphism
from $K\langle\langle Z,\theta_Z\rangle\rangle$ in $K\langle\langle
A,\theta\rangle\rangle$. Let $s$ be the morphism from $\langle
\beta_z(B)\rangle$ in $\M(Z,\theta_Z)$ defined by $s(zw)=z$ for each
$zw\in\beta_Z(B)$. We have
\[\begin{array}{rcccc}
\underline{\langle\beta_Z(B)\rangle}&=&s^{-1}(\underline{\M(Z,\theta_Z))}&=&
\sum_{w\in\M(Z,\theta_Z)}s^{-1}(w)\\
&=&\sum_{w\in\M(Z,\theta_Z)}\mu(w)&=&\mu(\underline{\M(Z,\theta_Z)})
\end{array}\]
Let $P(\theta_Z)$ be the polynomial such that
\[\underline{\M(Z,\theta_Z)}=\frac 1 {P(\theta_Z)}.\]
As $\mu$ is a continuous morphism, we have
\[\underline{\langle\beta_Z(B)\rangle}=\frac1{\mu(P(\theta_Z))}=\frac1{P
(\theta_{\beta_Z(B)})}\]
which is the characteristic series of $\M(\beta_Z(B),\theta_{\beta_Z(B)})$.
\CQFD

\begin{remark}
(i) Elimination in \cite{DK2} deals with the particular case when
$A-
B$ is totally non commutative. In this case $B$ is a TFSA of $A$.\\
(ii) As an example of other case, consider the independence alphabet due to
the graph\[\theta=a-b-c.\]
The monoid $\langle \beta_{a,b}(c)\rangle$ is free partially commutative,
its alphabet is $\beta_{a,b}(c)=\{b\}\cup\{ac^n/n\geq0\}$, its
independence graph is

\[\theta_{\beta_{a,b}(c)}=\begin{array}{ccccc}

& &ac& & \\
& &|&\ddots& \\
a&-&b&- &ac^n\\
& & & &\vdots
\end{array}\]
and
\[\underline{<\beta_Z(B)>}=\frac1{1-\left(b+\sum_{n\geq0}ac^n\right)+\sum_{n
\geq0}abc^n}.\]
\end{remark}

\section{Factorizations and bases of free partially commutative Lie
algebra}
\subsection{Transitive factorizations}
We recall some definitions given by Viennot in \cite{Vi1}.
\begin{defn}
Let $\M$ be a monoid, $\M'$ a submonoid of $\M$ and $\F=(\M_i)_{i\in J}$ a
factorization of $\M$. We denote $\F|_{\M'}=(\M_{i_k})_{k\in K}$ where
$K=\{k\in J|\M_k\subseteq \M'\}$ (in the general case it is not a
factorization).
\end{defn}
\begin{defn}\label{DDD5}
Let $\prec$ be the partial order on the set of all the factorizations of a
monoid
$\M$ defined by $\F=(\M_i)_{i\in J}\prec\F'=(\M'_i)_{i\in J'}$ ($\F'$ is
finer than $\F$) if and only
if $J'$ admits a decomposition $J'=\sum_{i\in J} J_i$ as an ordered sum of
intervals
such that for each $i\in J$, $(\M'_j)_{j\in J_i}$ is a factorization of
$\M_i$.
\end{defn}
The following property is straightforward.
\begin{prop}\label{PPP1}
Let $\F=(\M_i)_{i\in I}$ be a factorization and $\F'$ be a factorization
such
that $\F\preceq\F'$ then for each $i\in I$, $\F'|_{\M_i}$ is a
factorization of $\M_i$.
\end{prop}
\begin{defn}
Let $\B=(B_1,B_2)$ be a bisection and $\F=(Y_i)_{i\in J}$ a factorization.
We
say that $Y_i$ is {\rm cut} by $\B$ if and only if $\L_i(\B)=\langle
B_1\rangle\cap\langle Y_i\rangle$ and $\R_i(\B)=\langle B_2\rangle\cap
\langle Y_i\rangle$ are not trivial (i.e. not $\{1\}$).
\end{defn}

We need the following lemma.\\
\begin{lem}\label{LLL1}
Let $\B=(B_1,B_2)$ be a bisection of \hspace{1mm}$\M(A,\theta)$ and
$\F=(Y_i)_{i\in[1,n]}$ a
factorization with $n>1$, such that it exists a factorization
$\GG=(G_k)_{k\in K}$ with $\B,\F\preceq\GG$ then $\B\preceq\F$ if and
only if no $Y_i$ is cut by $\B$.
\end{lem}
{\bf Proof } We use the decomposition of K as an ordered sum of intervals
$K=J_1+J_2=\sum_{i\in[1,n]}I_i$
as in definition \ref{DDD5}. The assertion $(ii)$ implies
the existence of an integer $k\in [1,n]$ such that
$J_1=\sum_{i\in[1,k]}I_i$ and $J_2=\sum_{i\in[k+1,n]}I_i$. This allows us
to conclude.\CQFD
\begin{note}
In the preceding lemma, the existence of a common bound $\GG$ is essential
as shown by the 
following example (with $\M(A,\emptyset)=\{a,b,c\}^*$ and the rational
expressions written as in \cite{BR})
\[\B=(a,ba^*\cup ca^*)\mbox{ and } \F=(b,a,ab^+a^*\cup cb^*a^*)\]
No factor of $\F$ is cut by $\B$ and the two factorizations admit no common
upper bound.
\end{note}
In the sequel, we use the notion of a composition of factorizations as it is
defined
by
Viennot in \cite{Vi1}. We recall it here.
\begin{defn}
Let $\F=(\M_i)_{i\in I}$ be a factorization of a monoid $\M$ and for some
$k\in I$, $\F'=(\M'_i)_{i\in I'}$ a factorization of $\M_k$. The {\rm
composition} of $\F$ and $\F'$ is the factorization $\F'\circ
\F=(\M''_i)_{i\in I''}$ where $I''=I\cup I'-\{k\}$ is ordered by $i<j$ if
and only if
\begin{quote}
(i) ($i,j\in I$ and $i<_I j$) or ($i,j\in I'$ and $i<_{I'}j$))\\
(ii) $i\in I$, $i<_Ik$ and $j\in I'$ \\
(iii) $i\in I'$, $j>_Ik$ and $j\in I$
\end{quote}
and
\[\M''_i=\left\{\begin{array}{ll}
\M_i&\mbox{if }i\in I\\
\M'_i&\mbox{if }i\in I'\end{array}\right.\]
\end{defn}

\begin{defn}
A  {\rm transitive factorization} is a factorization which is composed
of transitive bisections (in finite number).
\end{defn}

\begin{lem}\label{LLL3}
Let $\F=(Y_i)_{i\in[1,p]}$ be a transitive factorization and
$\B=(B,\beta_Z(B))$ be a transitive bisection such that it exists a
factorization $\GG$ finer that $\B$ and $\F$. Then it exists at most one
$Y_i$ cut by $\B$ and for a such $i$ we have
\begin{quote}
(i) The subset $T=Y_i\cap\M(B,\theta_B)$ is a TFSA of $Y_i$ and
$\R_i(\B)$ is the right monoid of
the associated bisection.\\
(ii) The sequence $(Y_1, \dots ,Y_{i-1},T)$ is a transitive factorization
of
$\M(B,\theta_B)$.\\
(iii) The sequence $(\beta_{Y_i-T}(T),Y_{i+1}, \dots ,Y_p)$ is a transitive
factorization of $\M(\beta_Z(B),\theta_{\beta_Z(B)})$
\end{quote}
\end{lem}
{\bf Sketch of the proof} First it suffices to remark that, if $i>j$ are
two indices such that $Y_i$ and $Y_j$ are cut by $\B$ then
$\L_j(\B)\subseteq \M(B,\theta_B)\cap\M(\beta_Z(B),\theta_Z(B))=\{1\}$ and
this contradicts our hypothesis, hence $i=j$.\\
Let us prove assertion (i).\\
1) First, we remark that
\[\underline{\M(Y_i,\theta_{Y_i})}=\underline{\L_i(\B)}.\underline{\R_i(\B)}
\] and using the equality $\L_i(\B)=\M(T,\theta_T)$ we prove that
$\R_i(\B)=\langle \beta_{Y_i-T}(Y_i)\rangle$.\\
2) We show that if $T$ is not a TFSA of $Y_i$ then $B$ is not a TFSA of
$A$
and this implies the result.\\
Let us prove (ii) and (iii) by induction on $p$. If $p=1$ the result is
trivial. If $p>1$ , we can write $\F$ under the 
form
$\F=\F_1\circ\F_2\circ\B'$
where $\B'=(B',\beta_{Z'}(B'))$ is a transitive bisection,
$\F_1=(Y_1, \dots ,Y_{k})$ a transitive factorization of
$\M(B',\theta_{B'})$
and
$\F_2=(Y_{k+1}, \dots ,Y_p)$ a transitive factorization of
the monoid $\M(\beta_{Z'}(B'),\theta_{\beta_{Z'}(B')})$. If $\B=\B'$ the
result is trivial. If $\B\neq\B'$ , we have necessarily $B\subset B'$ or
$B'\subset B$. We suppose
that
$B'\subset B$ (the other case is symmetric), and we consider the transitive
trisection $(B',\beta_{B-B'}(B'),\beta_Z(B))$. Using the induction
hypothesis we find that
\[(Y_k, \dots ,Y_{i-1},T)\mbox{ and }(\beta_{Y_i-T}(T),Y_{i+1}, \dots
,Y_p)\]
are transitive factorizations
(respectively of the monoid $\M(\beta_{B-B'},\theta_{\beta_{B-B'}(B)})$ and
$\M(\beta_Z(B),\theta_{\beta_Z(B)})$). And then
\[(Y_1, \dots ,Y_{i-1},T)=\F_1\circ(Y_k, \dots
,Y_{i-1},T)\circ(B',\beta_{B-B'}(B))\]
is a transitive factorization. \CQFD
\begin{lem}\label{LLL4}
Let $\B=(B,\beta_Z(B))$ be a transitive bisection and
$\F=(Y_i)_{i\in[1,n]}$ be a transitive factorization such that
$\B\preceq\F$. Then the factorizations $\F|_{\M(B,\theta_B)}$ and
$\F|_{\M(\beta_Z(B),\theta_{\beta_Z(B)})}$ are transitive.
\end{lem}
{\bf Proof } We can prove the result by induction on $n$. \CQFD

\begin{prop}\label{PPP2}
Let $\F=(Y_i)_{i\in J}$ and $\F'=(Y'_j)_{j\in J'}$ be two finite transitive
factorizations such that it exists a factorization $\GG$ with
$\F,\F'\preceq\GG$ then it exists a transitive finite factorization $\GG'$
such that
\begin{quote}
(i) $\F,\F'\preceq \GG'\preceq\GG$\\
(ii) For each $j\in J$, $\GG'|_{\M(Y_j,\theta_{Y_j})}$ is a transitive
finite
factorization.\\
(iii) For each $j\in J'$, $\GG'|_{\M(Y'_j,\theta_{Y'_j})}$ is a transitive
finite
factorization.\\
\end{quote}

\end{prop}
{\bf Sketch of the proof } We set $J=[1,n]$, $J'=[1,n']$ and we prove the
result by induction on $n$. If $n=1$ the result is trivial. If $n=2$,
lemmas \ref{LLL1}, \ref{LLL3} and \ref{LLL4} give us the proof. If $n\geq
2$, we set $\F=\F_1\circ\F_2\circ\B$ where $\B=(B,\beta_Z(B))$ is a
transitive bisection of $\M(A,\theta)$, $\F_1$ a transitive factorization
of
$\M(B,\theta_B)$ and $\F_2$ a transitive factorization of
$\M(\beta_Z(B),\theta_{\beta_Z(B)})$. Using lemmas \ref{LLL1}, \ref{LLL3}
and
\ref{LLL4} we define a factorization
\[\F''=\left\{\begin{array}{ll}
\F'&\mbox{If } B\preceq\F'\\
(Y'_{1}, \dots ,Y'_{i-1},T,\beta_Z(T),Y'_{i+1}, \dots
,Y'_{n'})&\mbox{Otherwise}
\end{array}
\right.\]
such that $\F',\B\preceq\F''\preceq\GG$, $\F''|_{\M(Y'_j,\theta_{Y'_j})}$
is
transitive for each $j\in [1,n]$ (in fact this factorization is trivial for
all $j\in[1,n]$ but at most one for which it is a transitive bisection),
$\F''|_{\M(B,\theta_B)}$ and $\F''|_{\M(\beta_Z(B),\theta_{\beta_Z(B)})}$
are
transitive. Using the induction hypothesis we can construct two
factorizations $\F''_1$ and $\F''_2$ such that
\[\F_1,\F''|_{\M(B,\theta_B)}\preceq\F''_2\preceq\GG|_{\M(B,\theta_B)}\]
and
\[\F_2,\F''|_{\M(\beta_Z(B),\theta_{\beta_Z(B)})}\preceq\F''_2\preceq\GG|_{\
M(\beta_Z(B),\theta_{\beta_Z(B)})}\] and satisfying (ii) and (iii). We set
$\GG'=\F''_1\circ\F''_2\circ\B$, then $\F,\F'\preceq\GG'\preceq\GG$ and the
induction hypothesis, the construction of $\F''$ and lemma \ref{LLL4}
allow us to conclude. \CQFD
\begin{cor}\label{CCC1}
Let $\F=(Y_i)_{i\in I}\preceq\F'$ be two transitive finite factorizations
then for each $i\in I,\F'|_{\M(Y_i,I_{Y_i})}$ is a transitive finite
factorization.
\end{cor}
{\bf Proof } It suffices to use proposition \ref{PPP2} with
$\F,\F'\preceq\F'$.\CQFD

The following definition is an adaptation to partial commutations of a
definition
given by Viennot in \cite{Vi1}.
\begin{defn}\label{D1}
A factorization $(Y_i)_{i\in I}$ of $\M(A,\theta)$ has {\rm locally the
property $\goth P$}
if and only if for each finite subalphabet $B\subset A$ and $n\geq 0$ it
exists a factorization $(Y'_i)_{i\in I'}$ with the property $\goth P$ such
that there is an strictly increasing mapping $\phi: I'\rightarrow I$
satisfying
\begin{center}
$\displaystyle Y'_i\cap B^{\leq n}=
Y_{\phi(i)}\cap B^{\leq n}$
and $ Y_j\cap B^{\leq n}=\emptyset$ if $j\not\in \phi(I')$
\end{center}
\end{defn}

\begin{defn}
We denote $CLTF(A,\theta)$ the set
of the complete locally transitive finite factorizations.
\end{defn}

\begin{example}
Consider the following independence graph
\[a-b-c-d.\]
We construct a complete locally transitive finite factorization $\F$ as
follow.
We eliminate successively the traces $c, ac^2, b, d, ac$ and $a$. So we
have \[M(A,\theta)=c^*.(ac^2)^*.b^*.d^*.(ac)^*.a^*.M\] where $M$ is a
(non-commutative) free monoid. It suffices to take a Lazard factorization
on $M$ to construct a complete locally transitive finite factorization of
$\M(A,\theta)$.\\
We can remark that one can not obtain this
factorization using only transitive bisections with a non commutative right
member. Examining all the transitive bisections of this kind
$$
\begin{array}{ll}
1.\quad\B_1=(\{a,c\},\beta_{b,d}(a,c))&5.\quad\B_5=(\{a,c,d\},\beta_b(a,c,d)
)\\
2.\quad\B_2=(\{b,c\},\beta_{a,d}(b,c))&6.\quad\B_6=(\{b,c,d\},\beta_a(b,c,d)
)\\
3.\quad\B_3=(\{b,d\},\beta_{a,c}(b,d))&7.\quad\B_7=(\{a,b,d\},\beta_c(a,b,d)
),\\
4.\quad\B_4=(\{a,b,c\},\beta_d(a,b,c))&
\end{array}$$
we can easily prove that $\F$ could not be written like
$\F=\F_1\circ\F_2\circ\B_i$ with 
$i\in \{1,2,\dots,7\}$.
\end{example}

\subsection{Transitive elimination in $L_K(A,\theta)$}

The algebra of  trace polynomials $K\langle A,\theta\rangle=K[\M(A,\theta)]$
endowed with the 
classical Lie bracket is a Lie algebra (\cite{Bo,DK1,Lo} when
$\theta=\emptyset$). The free partially commutative Lie algebra is at first
defined as the free object with respect to the given commutations
\cite{DK1}. Here we will use 
its realisation as the sub-Lie algebra of $K\langle A,\theta\rangle$
generated by the letters 
\cite{La2}. We will denote it $L_K(A,\theta)$. One can show that this
definition is equivalent to 
$L_K(A,\theta)=L_K(A)/_{I_\theta}$ where $L_K(A)$ is the free Lie algebra
and $I_\theta$ is the 
Lie ideal of $L_K(A)$ generated by the polynomials $[a,b]$ with $(a,b)\in
\theta$.
The following theorem proves that elimination in $L_K(A,\theta)$ and
transitive factorization of $\M(A,\theta)$ occur under the same condition.
\begin{thm}\label{LieTh1}
Let $(B,Z)$ be a partition of $A$
\begin{quote}
(i) We have the decomposition
\[L_K(A,\theta)=L_K(B,\theta_B)\oplus J\]
where $J$ is the Lie ideal generated (as a Lie algebra) by
\[\tau_Z(B)=\{[ \dots [z,b_1], \dots b_n]\quad|\quad zb_1 \dots b_n\in
\beta_Z(B)\}. \]
(ii) The subalgebra $J$ is a free partially commutative Lie algebra if
$B$ is a TFSA of $A$.\\
(iii) Conversely if $J$ is a free partially commutative Lie algebra with
code $\tau_Z(B)$ then 
$B$ is a TFSA.
\end{quote}
\end{thm}
{\bf Proof} $(i)$ We have the classical Lazard bisection
\[L_K(A)=L_K(B)\oplus L_K(T_Z(B))\]
where $T_Z(B)=\{[ \dots [z,b_1], \dots ],b-n]\quad |\quad z\in Z, b_1,
\dots ,b_n\in B\}$. Then using the
natural mapping $L_K(A)\rightarrow L_K(A,\theta)$ (as
$[ \dots [z,b_1], \dots ],b_n]$ maps to $0$ if $zb_1 \dots
b_n\not\in\beta_Z(B)$) we
get
the claim.\\
$(ii)$ The proof goes as in \cite{DK2}, due to the fact that, for a TFSA,
(c) below still holds, we 
sketch the proof.\\
Define a
mapping $\partial_b$ from $\beta_Z(B)$ to
$L_K(\beta_Z(B),\theta_{\beta_Z(B)})$
by
\[\partial_b=\left\{
\begin{array}{ll}
zwb&\mbox{if } zwb\in \beta_Z(B),\\
0&\mbox{otherwise}.
\end{array}\right. \]
a) We prove that if $B$ is TFSA, $\partial_b$ can be extended as a
derivation of the Lie algebra $L_K(\beta_Z(B),\theta_{\beta_Z(B)})$.\\
b) We define $\partial$ a mapping from $B$ to
$Der(L_K(\beta_Z(B),\theta_{\beta_Z(B)}))$ by $\partial(b)=\partial_b$ and
we extend it as a Lie morphism from $L_K(B,\theta_B)$ into
$Der(L_K(\beta_Z(B),\theta_{\beta_Z(B)}))$.\\
c) We prove that the semi-direct product $L_K(B,\theta_B)\propto_\partial
L_K(\beta_Z(B),\theta_{\beta_Z(B)})$ and the Lie algebra $L_K(A,\theta)$
are
isomorphic using the universal property of the latter. Hence, $J$ is a free
partially commutative 
Lie algebra
isomorphic to $L_K(\beta_Z(B),\theta_{\beta_Z(B)})$.\\
(iii) If the dependence graph admits the following subgraph
\[z-b_1- \dots -b_n-z'\]
with $b_i\in B$ and $z,z'\in Z$ we have the identity
\[[z,[[ \dots [z',b_n] \dots ,b_2],b_1]]=[ \dots [z',b_n] \dots
b_2],[z,b_1]].\]
Which implies that $\tau_Z(B)$ is not a code for $J$.\CQFD

\subsection{Construction of bases of $L_K(A,\theta)$}
In this section, we define a class of bases which contains the bases found
by Duchamp and Krob in \cite{DK1}, \cite{DK2}
and \cite{DR} using chromatic partitions and the partially
commutative Lyndon bases found by Lalonde (see Lalonde \cite{La},
Krob and Lalonde \cite{KL}).\\
\begin{defn}
Let $\F=(Y_i)_{i\in[1,n+1]}$ be a finite transitive factorization. We
denote
$\widetilde \F$, the set of the n-uplets $(\B_1, \dots ,\B_n)$ of
transitive
bisections such that $\F=\B_n\circ \dots \circ\B_1$.\\
Let $\F$ be a transitive factorization and ${\goth
f}=(\B_1, \dots ,\B_n)\in\widetilde\F$, we denote ${\goth f}\B_n^{-
1}=(\B_1, \dots ,\B_{n-1})$.
\end{defn}
In the following, if $\F$ is the sequence $(Y_i)_{i\in J}$, we denote
$\Cont\F=\bigcup_{i\in 
J}Y_i$ as in \cite{Vi1}.
\begin{defn}
Let $\F=(Y_i)_{i\in[1,n+1]}$ be a finite transitive factorization. A
bracketing of $\F$ along ${\goth f}\in\widetilde\F$ is a mapping
$\Pi_{\goth f}$ from $\Cont\F$ to $L_K(A,\theta)$
inductively defined as follows.
If $n=1$, then $\goth f$ is a sequence of length 1 under the form
$((B,\beta_Z(B)))$
and
\[\Pi_{\goth f}(w)=\left\{\begin{array}{l}
w \mbox{ if } w\in B,\\
{[ \dots [z,b_1] \dots b_k]} \mbox{ if } w=zb_1 \dots b_k \in
\beta_Z(B)\mbox{ and }
z\in Z.
\end{array}\right. \]
If $n>1$, let ${\goth f}=(\B_1, \dots ,\B_n)\in\widetilde\F$.
We set
$\B_{n-1}\circ \cdots \circ \B_1=(Y'_i)_{i\in[1,n]}$ and $j\in[1,n]$ such
that $\B_n=(Y''_j,\beta_{Y'_j-Y''_j}(Y''_j))$ (remark that in this case, one
has $\Cont\F-
\Cont\B_{n-1}\circ\cdots\circ\B_1=\beta_{Y'_j-Y''_j}(Y''_j)$).
And
\[\Pi_{\goth f}(w)=\left\{\begin{array}{ll}
\Pi_{{\goth f}\B_n^{-1}}(w) \mbox{ if } w\in
\Cont\B_{n-1}\circ\cdots\circ\B_1,&\\
{[ \dots [\Pi_{{\goth f}\B_n^{-1}}(y_1),\Pi_{{\goth f}\B_n^{-1}}
(v_1)], \dots \Pi_{{\goth f}\B_n^{-1}}(v_k)]}&\mbox{ if }
w=y_1v_1 \dots v_k,\\
&w \in \beta_{Y'_j-Y''_j}(Y''_j),\\
& y_1\in Y'_j-Y''_j\\
&\mbox{and } v_1, \dots, v_k\in Y''_j.
\end{array}\right. \]
\end{defn}
Using theorem \ref{LieTh1} in an induction on $n$ we prove the following
proposition.
\begin{prop}\label{LieP1}
Let $\F=(Y_i)_{i\in[1,n]}$ be a transitive factorization. For each ${\goth
f}\in\widetilde\F$, we have the following decomposition
\[L_K(A,\theta)=\displaystyle\bigoplus_{i\in[1,n-1]}L_K(\Pi_{\goth
f}(Y_i),\theta_i)\]
where
\[\theta_i=\{(\Pi_{\goth f}(y_1),\Pi_{\goth f}(y_2))|(y_1,y_2)\in
\theta_\M\}. \]
\end{prop}
\begin{defn}
Let $\F=(Y_i)_{i\in J}$ be a locally transitive finite factorization, a
{\it bracketing} of $\F$ is
a mapping
$\Pi$ from $\bigcup_{i\in J}Y_i$ to $L_K(A,\theta)$
such that for each finite subalphabet $B\subset A$ and each integer $n\geq
0$, it exists a transitive finite factorization $\F_{n,B}=(Y_i^{n,B})_{i\in
J_{n,B}}$ and ${\goth f}_{n,B}\in\widetilde\F_{n,B}$ such that for each
$t\in \Cont\F_{n,B}\cap B^{\leq n}$,
$\Pi(t)=\Pi_{{\goth f}_{n,B}}(t)$.
\end{defn}

\begin{lem}\label{LF3}
Let $\F=(Y_i)_{i\in J}\preceq \F'$ be two finite transitive factorizations.
Then, for each ${\goth f}\in \widetilde\F$, it exists ${\goth
f'}\in\widetilde\F'$ such that for each $t\in \Cont\F$,
$\Pi_{\goth f}(t)=\Pi_{\goth f'}(t)$.

\end{lem}
{\bf Proof } It is a direct consequence of corollary \ref{CCC1}. \CQFD
\begin{thm}
Let $(A,\theta)$ be an independence alphabet. Each locally finite
transitive
factorization of $\M(A,\theta)$ admits a bracketing.
\end{thm}
{\bf Proof } Let $\F=(Y_i)_{i\in J}$ be a locally finite transitive
factorization. Using 
proposition \ref{PPP2}, one can construct a sequence of finite transitive
factorizations 
$(F_{n,B})_{n\in\N, B\subset A\atop Card B<\infty}$ such that
\begin{enumerate}
\item if $n\leq n'$ and $B\subset B'$ then
\[\F_{n,B}\preceq\F_{n',B'},\]
\item for each $n\geq 0$ and $B\subset A$
\[\F_{n,B}\preceq\F,\]
\item for each $n\geq 0$ and each finite subalphabet $B$, if we set 
$\F_{n,B}=(Y^{n,B})_{i\in[1,k_{n,B}]}$, it exists a strictly increasing
mapping $\phi_{n,B}$ from 
$[1,k_{n,B}]$ to $J$ verifying
\[\M(Y_i^{n,B},\theta_{Y_i^{n,B}})\cap B^{\leq
n}=\M(Y_{\phi_{n,B}(i)},\theta_{Y_{\phi_{n,B}(i)}})\]
and
\[j\notin\phi_{n,B}([1,k_{b,B}])\Rightarrow \M(Y_j,\theta_{Y_j})\cap B^{\leq
n}=\{1\}.\]
\end{enumerate}
By lemma \ref{LF3}, we can define for each $n>0$ and each finite subalphabet
$B$ of $A$ a 
sequence ${\goth f}_{n,B}\in\tilde\F_{n,B}$ such that for each $m<n$,
$B'\subset B$ and $t\in 
\Cont\F_{m,B'}\cap B'^{\leq n}$ we have $\Pi_{{\goth
f}_{m,B'}}t=\Pi_{{\goth f}_{n,B}}t$.\\
Thus, we can define $\Pi$ as the mapping from $\Cont\F$ into $L_K(A,\theta)$
such that $\Pi 
t=\Pi_{{\goth f}_{|t|,\Alph(t)}}t$.\CQFD \\
We have easily the following result.
\begin{prop}
Let $\F=(\{l_i\})_{i\in I}\in CLTF(A,\theta)$ and $\Pi$ be a bracketing of
$\F$ then the family
$(\Pi(l_i))_{i\in I}$ is a basis of $L_K(A,\theta)$ as K-module.
\end{prop}
\begin{example}
We set $A=\{a,b,c,d\}$ and $\theta=a-b-c-d$. We construct locally (for
$n\leq 3$) the following basis.
\end{example}
[[a,d],b], [[a,d],d], [[a,d],a], [a,d], [a,[a,c]], a, [a,c], [[a,c],c],
[[a,d],c], [b,d], [[b,d[,b[, [[b,d],d], b, c, d.

\section{The case of the group}

The free partially commutative group \cite{DT} can be defined by the
presentation
\[\F(A,\theta)=<A;\{ab=ba\}_{(a,b)\in\theta}>_{gr}.\]
To extend the elimination process, we need the alphabet of the inverse
letters. Recall that one 
can construct the free partially commutative group using "reduced" traces
\cite{DR,DK3}. If $A$ 
is an alphabet, we define $\tilde A=A\cup\OV A$ where $\OV A=\{\OV a\}_{a\in
A}$ is a disjoint 
copy of $A$. The set $\tilde A$ is then provided with the involution
$x\rightarrow \OV x$ such that 
$\OV{\OV x}=x$. Thus, $\theta$ is extended by
\[\tilde\theta=\{(x,y)\in \tilde A^2|\{(x,y),(\OV x,y),(x,\OV y),(\OV x,\OV
y)\}\cap\theta\neq\emptyset\}.\]
We define the natural mapping $s_0:\tilde A\rightarrow\F(A,\theta)$ such
that $s_0(a)=a$ and 
$s_0(\OV a)=a^{-1}$ for each letter $a\in A$.\\
As $s_0$ is compatible wih the commutations of $\tilde\theta$ (i.e. if
$(x,y)\in\tilde\theta$ 
then $s_0(x)s_0(y)=s_0(y)s_0(x)$), one has the factorization
\[
\begin{array}{c@{\hskip 1cm}c}
\rnode{a}{\tilde A}&\rnode{c}{\F(A,\theta)}\\[1cm]
\rnode{b}{\M(\tilde A,\tilde\theta)}&\hfill .
\end{array}
\ncline{->}{a}{c}\Aput{s_0}
\ncline{<-}{c}{b}\Aput{s}
\ncline{->}{a}{b}
\]
The mapping $s$ is onto. For each $g\in\F(A,\theta)$, it exists an unique
preimage with minimal 
length in $s^{-1}(g)$, this element is called "reduced expression" of $g$
(the subset of these 
traces will be denoted by $red(\tilde A,\tilde\theta)$). The links with the
bisections 
$(B,\beta_Z(B))$ is given by the following.
\begin{lem}
Let $B\subset A$, $z\not\in B$ and $w\in red(\tilde B,\tilde\theta_{\tilde
B})$. The following 
assertions are equivalent.
\begin{enumerate}
\item $\OV wzw$ is a reduced trace (i.e. $\OV wzw\in red(\tilde
A,\tilde\theta)$).
\item $zw\in\beta_z(\tilde B)$.
\end{enumerate}
\end{lem}
{\bf Proof} Straightforward, using the criterion given in \cite{DK3}:\\ Let
$t=a_1a_2\dots 
a_n\in\M(\tilde A,\tilde\theta)$, $t$ is not a reduced trace if and only if
it exists $1\leq 
i<j\leq n$ with $a_i=\OV a_j$ and such that for each $k$, $i<k<j$, 
$(a_k,a_i)\in\tilde\theta$.\cqfd\\ \\
We denote $\beta_Z^R(\tilde B)$ the set $\beta_Z(\tilde B)\cap red(\tilde
A,\tilde\theta)$ with the commutation 
$\tilde\theta_{\beta_Z^R(\tilde B)}$ provided by the definition
\ref{ComThe}. One 
has an analogue of the theorem \ref{LieTh1}.
\begin{prop}
Let $(B,Z)$ be a partition of $A$.
\begin{enumerate}
\item[(i)] One has the decomposition as the semi direct 
product$$\F(A,\theta)=\F(B,\theta_B)\ltimes H_Z $$
where $H_Z$ is the normal subgroup generated by $Z$. It is the subgroup
generated by
$$\rho_Z(B)=\{w^{-1}zw|zw\in\beta_z^R(\tilde B)\}$$
\item[(ii)]\label{A2} The subgroup $H_Z$ is free partially commutative for
the code $\rho_Z(B)$ and the 
commutations $$\hat\theta_\rho:=\{(t,t')\in\rho_Z(B)^2|tt'=t't\mbox{ and }
t\neq t'\}.$$
\item[(iii)]The natural mapping $\alpha:\F(\beta_Z^R(\tilde
B),\tilde\theta_{\beta_Z^R(\tilde 
B)})\rightarrow H_Z$ is one to one if and only if $B$ is TFSA.
\end{enumerate}
\end{prop}
{\bf Proof} (i) The decomposition given by (i) is the image of the non
commutative Lazard elimination in the free group. The 
unicity of the decomposition with respect to the semidirect product can be
obtain (as in the 
classical case) by sending all the element of $Z$ to one.\\
(ii) Let $\hat\rho=\{a_t\}_{t\in\rho_Z(B)}$ be an alphabet and $\hat\theta$
be the commutation 
relation defined by $(a_t,a_{t'})\in\hat\theta$ if and only if $t\neq t'$
and $tt'=t't$.\\
For each $b\in B$, we define the mapping $\sigma_b:\hat\rho \rightarrow
\hat\rho$ by 
$\sigma_b(a_t)=a_{b^{-1}tb}$. Remarking that $b^{-1}tb$ belongs to
$\rho_Z(B)$ and then 
$(a_t,a_{t'})\in\hat\theta$ implies
$(\sigma_b(a_t),\sigma_b(a_{t'}))\in\hat\theta$, this mapping 
can be extended in an automorphism $\sigma_b$ of $\F(\hat\rho,\hat\theta)$.
Let $\sigma$ be the 
mapping from $B$ to $Aut(\F(\hat\rho,\hat\theta))$ defined by
$\sigma(b)=\sigma_b$. As 
$\sigma_b\sigma_{b'}=\sigma_{b'}\sigma_b$ when $(b,b')\in\theta_B$, $\sigma$
can be extended as a 
morphism from $\F(B,\theta_B)$ in $Aut(\F(\hat\rho,\hat\theta))$. Using the
same proof than in 
theorem \ref{LieTh1}, we find that the semidirect product
$\F(B,\theta_B)\propto_\sigma 
\F(\hat\rho,\hat\theta)$ and $\F(A,\theta)$ are isomorphic.\\ \\

(iii) Suppose that $B$ is not TFSA then it exists a $(z_1,z_2)\in\theta_Z$
and a minimal path in 
the non commutation graph \[z_1-a_1-\cdots-a_k-c-b_l-\cdots-b_1-z_2.\]
Let $r_1=z_1a_1\cdots a_k$ and $r_2=z_2b_1\cdots b_l$. Due to the fact that
the chain is of minimal length, one has 
$(r_1,r_2)\in\tilde\theta_{\beta_Z^R(\tilde B)}$ and
$\alpha(r_1)\alpha(r_2)=\alpha(r_2)\alpha(r_1)$. But $r_1c$ and $r_2c$ 
do not commute and their images $\alpha(r_1c)=c^{-1}\alpha(r_1)c$ and
$\alpha(r_2c)=c^{-1}\alpha(r_1)c$ do. This proves 
that $\alpha$ is not one to one. \\
The converse follows from the fact that, when $B$ is a TFSA, the commutation
graph 
$(\beta_Z^R(\tilde B),\tilde\theta_{\beta_Z^R(\tilde B)})$ and
$(\rho_Z(B),\hat\theta_\rho)$ are 
obviously isomorphic.
\cqfd
\begin{note}
In general $\alpha$ is into.
\end{note}

\end{document}